\input amstex
\input xy
\xyoption{all}
\input epsf
\documentstyle{amsppt}
\document
\magnification=1200
\NoBlackBoxes
\nologo
\hoffset=1.5cm
\voffset=1truein

\pageheight{16cm}


\bigskip

 \centerline{\bf NUMBERS AS FUNCTIONS\footnotemark1}
 \footnotetext{Based on talks at the {\it International Workshop on $p$--adic methods for modelling of
complex systems}, Bielefeld, April 15--19, 2013, and 
 at  {\it Journ\'ees Arithm\'etiques}, Grenoble, June 2--5, 2013.} 

\medskip

\centerline{\bf Yuri I. Manin}

\smallskip

\centerline{Max--Planck--Institut f\"ur Mathematik, Bonn, Germany}
  
  \bigskip
  
 {\it ABSTRACT.  In this survey I discuss A.~Buium's theory of ``differential equations in the $p$--adic direction'' ([Bu05])
 and its interrelations  with  ``geometry over  fields with one element'', on the background
 of various approaches to $p$--adic models in theoretical physics (cf. [VlVoZe94], [ACG13]).} 
  
  \bigskip
  
 \centerline{\bf Introduction}
 
 \medskip 

One of the most beautiful (arguably, {\it the} most beautiful) mathematical formulas is
Euler's identity
$$
e^{\pi i}=-1 .
\eqno(0.1)
$$
It connects four numbers  $\pi=3,1415912\dots$, $e=2,71828\dots$, $i=\sqrt {-1}$, and $-1$ itself,
and has a very strong physical flavor being the base of the universal  principle
of ``interference of probability amplitudes'' in quantum mechanics and quantum field theory.
The ``$-1$'' in the right hand side of (0.1) shows how two quantum states
with opposite phases may annihilate each other after superposition.

\smallskip

On the other hand, of these four numbers $\pi, e, i, -1$, {\it only} $\pi$ looks as something
similar to a ``physical constant'' in the sense that it can be (and was) measured, with a certain
approximation. 

\smallskip
Moreover, the traditional names of the respective classes of numbers,
which we nowadays tend to perceive as mathematical terms introduced by precise definitions in courses of calculus, --
{\it irrational, transcendent, imaginary, negative,} -- in the course of history conveyed
the primeval bafflement of the rational mind, discovering these numbers but reluctant to accept them.

\smallskip

We may recall that at the time of their discovery these numbers had very different sources of justification:
$\pi$ in Euclidean geometry (which describes essentially kinematics of solids in gravitational vacuum), $-1$
in commerce (``debt''), $e$ in the early history of computer science (Napier's implementation of the discovery that
 a specific precomputation can facilitate everyday tasks of multiplication), $i$ in the early history of polynomial
 equations.

\smallskip

When I was asked to deliver a talk at the  Workshop on $p$--adic Methods for Modelling of
Complex Systems, I decided to present first a  $p$--adic environment of $\pi$ and $e$.

\smallskip

Probably, the earliest ``arithmetic'' formula involving $\pi$ is due to Euler (as well as (0.1)):
$$
\frac{\pi^2}{6} =\prod_p (1-p^{-2})^{-1}.
\eqno(0.2)
$$
However, it involves {\it all} primes $p$ simultaneously, and in fact, can be best understood
as a fact from {\it ad\`elic} geometry. As such, it looks as a generalisation of the
simple--minded product formula $\prod_v |a|_v=1$ valid for all $a\in \bold{Q}^*$,
where $v$ runs over all valuations of $\bold{Q}$, $p$--adic ones and archimedean one.
To be more precise, (0.2) expresses the fact that the natural adelic measure of
$SL(2,A_{\bold{Q}})/SL(2,\bold{Q})$ equals 1. For some more details, cf. [Ma89],
where it was suggested that fundamental quantum physics might be related to number theory
via this ad\`elic philosophy, {\it ``democracy of all valuations''}, and the exclusive use of real and complex numbers
in our standard formalisms  is the matter of tradition, which we now try to overcome by
replacing ``the first among equals''  archimedean valuation by an arbitrary non--archimedean one.

\smallskip

Now we turn to $e$. Here, as the discoverer of $p$--adic numbers Kurt Hensel himself remarked,
we have a candidate for $e^p$ in each $p$--adic field, since the (archimedean) series for
$e^p$ converges also $p$--adically:
$$
e^p=\sum_{n=0}^{\infty} \frac{p^n}{n!}.
\eqno(0.3)
$$
Since the root of degree $p$ of the right hand side  of (0.3) understood as
$p$--adic number generates an extension of $\bold{Q}_p$ of degree $p$,
there can be no algebraic number with such local components.

\smallskip

This argument looks tantalisingly close to a proof of transcendence of $e$,
although, of course, it is not one. On the other hand, I do not know any ad\`elic formula
involving $e$ in such a way as (0.2) involves $\pi$.

\smallskip

In this survey, I proceed with discussion consisting of three main parts.

\medskip

 {\it A. I will describe a class of numbers (including transcendental ones) 
relevant for Quantum Field Theory in the sense that they define
the coefficients of perturbative series for Feynmann's path integrals.
These numbers are called {\it (numerical) periods}, they were introduced and studied in 
[KoZa01].}

\smallskip
 Roughly speaking, numerical periods are values at algebraic points  of 
certain {\it multi--valued transcendental functions}, naturally defined on various moduli spaces,
and also traditionally called (functions--)periods.

\smallskip

These functions--periods satisfy differential equations of Picard--Fuchs type,
and such equations furnish main tools for studying them.

\smallskip

In the second part of this survey, I focus on the following program:

\medskip

{\it B. For a prime $p$, numerical periods also can be considered as solutions
of  ``differential equations in the $p$--adic direction''}. 

\smallskip
The whole machinery
of such differential equations was suggested and developed by Alexandru Buium,
cf. his monograph [Bu05], and I briefly review it. I use the catchword ``numbers as functions'' to name this analogy.
\smallskip

Alexandru Buium
has convincingly shown that  the right analog 
of the $p$--adic derivation is (a natural generalization of) the Fermat quotient  $\delta_p(a):= (a-a^p)/p$
initially defined for $a\in \bold{Z}$. 
Unexpectedly, this formal idea had rich consequences: Buium was able to construct analogs of classical jet spaces
``in the $p$--adic direction'',  together with a theory of functions on these jet spaces,
containing an incredible amount of analogs of classical constructions traditionally
requiring calculus. 

\smallskip
Those numerical periods that were already treated by Buium include
periods of abelian varieties defined over number (or even $p$--adic) fields. (But the reader
should be aware that, in the absence of uniformization, this last statement only very crudely describes a
pretty complicated picture;  see more details in the main text.)

\medskip

{\it C. For Buium's differential equations, ``constants in the $p$--adic direction'' turn out to be  roots of unity and zero:
Teichm\"uller's representatives of residue classes modulo $p$.}

\smallskip

Until recently, algebraic geometry over such constants was
motivated by very different insights: for a more detailed survey cf. [Ma95], [Ma08]. It is known as ``theory of the field $\bold{F}_1$''.

\smallskip
Briefly, this last field of inquiry is focused on the following goal: to make
the analogy between, say, $Spec\,\bold{Z}$ (or spectra of rings of algebraic integers) on the one hand,
and algebraic curves over finite fields, on the other hand, so elaborated
and precise that one could use a version of the technique of Andr\'e Weil, Alexander Grothendieck and Pierre Deligne
in order to approach Riemann's conjecture for Riemann's zeta and similar arithmetic functions.
\smallskip

The solid bridge between $\bold{F}_1$--geometry and arithmetic differential equations
was constructed by James Borger: cf. [Bor11a,b], [Bor09], [BorBu09].
Roughly speaking, in order to define the $p$--adic derivative $\delta_p$ of elements
of a commutative ring $A$, one needs {\it a lift} of the Frobenius map, that is an endomorphism $a\mapsto F(a)$,
such that $F(a) \equiv a^p\,\roman{mod}\,p.$ Borger remarked that a very natural system
of such lifts for all $p$ simultaneously is encoded in the so called {\it psi--structure}
or its slight  modification, {\it lambda--structure,} and then suggested to consider such a structure
as {\it descent data on $\roman{Spec}\,A$ to $\bold{F}_1$}.  
A related notion of ``cyclotomic coordinates'' in $\bold{F}_1$was independently suggested in [Ma08].
In particular, $a\in A$ is a cyclotomic co--ordinate (wrt a prime $p$) if $F(a)=a^p.$
I will return to these ideas in the last part of this survey.

\smallskip

Finally, I should mention that there exists a very well developed deep theory
of  ``$p$--adic periods'' for algebraic varieties defined over $p$--adic fields that
replaced the classic integration of differential forms over topological cycles
with comparison of algebraic de Rham and \'etale  cohomology theories:
see [Fa88] and a recent contribution and brief survey [Be11]. Periods in this setting
belong to a very big Fontaine's field $B_{dR}$. 
The approach to  periods via Buium's $p$--adic geometry that we describe in this survey has a very different flavour.
It would certainly be important to find connections between the two theories.

\bigskip

\centerline{\bf 1. Periods}

\medskip

{\bf 1.1. Numerical periods.}  M. Kontsevich and D.  Zagier introduced  an important
subring  $\Cal{P} \subset\bold{C}$ containing all algebraic numbers and a lot of numbers important in physics 
(see [KoZa01]).

\medskip

{\bf 1.1.1. Definition.} {\it $\alpha \in \Cal{P}$ if and only if the
real and imaginary parts of  $\alpha$  are values of
absolutely convergent  integrals
of functions in $\bold{Q}(x_1,\dots , x_n)$ over chains in
$\bold{R}^n$ given by polynomial (in)equalities with coefficients in $\bold{Q}$.}

\medskip

{\bf 1.1.2. Examples.} a) All algebraic numbers are periods.

\smallskip
b)  $\pi =\int\int_{x^2+y^2\le 1} dxdy$.

\smallskip

c)  $\Gamma\,(p/q)^q \in\Cal{P}$.
\medskip

It is not difficult to prove that periods form a subring of $\bold{C}$.
Feynman integrals (of a certain class) are periods.
But it is still not known whether $\pi^{-1}$, $e$, or Euler's constant
$\gamma$  are periods (probably, not). There is a close connection between periods and Grothendieck motives
(see [KoZa01]), and $2\pi i$ corresponds to the  Tate's motive. Since
in the motivic formalism one formally inverts the Tate motive, it is also
useful to extend the period ring by  $(2\pi i)^{-1}$.

\medskip 

d) The multiple $\zeta$--values  (Euler) 
$$
\zeta(n_1, ..., n_m) = \sum_{0< k_1 < ... < k_m}
\frac{1}{k_1^{n_1} ...\, k_m^{n_m} }\, , \qquad n_i\ge 1, n_m >1 \,.
\eqno(1.1)
$$
are periods.

\smallskip

In order to see it, we reproduce the Leibniz and Kontsevich integral formula for them.
\smallskip
Let  $n_1,\dots ,n_m$  be positive integers as in (1.1). Put
$n:=n_1+\dots +n_m$, and  $\underline \varepsilon:=
(\varepsilon_1, ..., \varepsilon_n)$ where $\varepsilon_i = 0$ or $1$,
and
$\varepsilon_i=1$ precisely when
$i\in\{1,n_1+1,n_1+n_2+1,\dots ,n_1+\dots + n_{m-1}+1\}$.
Furthermore, put
$$
\omega(\underline \varepsilon):= \frac{dt_1}{t_1-\varepsilon_1}\wedge
...
\wedge \frac{dt_n}{t_n-\varepsilon_n}
$$
and
$$
\Delta^0_n:= \{(t_1^0, \dots, t_n^0) \in \bold{R}^n\quad | \quad 0 < t_1^0 < \dots <
t_n^0 < 1\}
$$
Then we have
$$
\zeta (n_1,\dots ,n_m)= \zeta (\underline \varepsilon )
= (-1)^m\int_{\Delta^0_n} \omega (\underline \varepsilon ).
$$

For further details, see [GoMa04], where the mixed motives associated with these periods were identified:
they are constructed using moduli spaces $\overline{M}_{0,n}$ and their
canonical stratifications.

\medskip
 
 {\bf 1.2. Periods--functions.} Sometimes we may introduce parameters in the 
 description of elements of $\Cal{P}$ sketched above and thus pass to the study of
 periods as functions. To this end, it is first convenient
 to rewrite the definition in  a more formal algebraic--geometric framework
 as was already done in [KoZa01], sec. 4.1.
 
 \smallskip
 
 Consider a quadruple $(V, D, \omega , \gamma )$. Here $V$ is a smooth
 algebraic variety of pure dimension $n$, endowed with divisor $D$ with normal
 crossings, $n$--form $\omega$ regular outside $D$, and 
a homology class $\gamma \in H_n(V(\bold{C}), D(\bold{C});\bold{Q})$.
Moreover, $(V,D,\omega )$ must be defined over $\bold{Q}$, and the integral
$\int_{\gamma}\omega$ must converge. Then the set of
such integrals coincides with the period ring $\Cal{P}$ defined above.

\smallskip

It is now clear how to relativise this definition, replacing $V$ by
a relatively smooth morphism  $f:\,\Cal{V}\to S$ defined over $\bold{Q}$,
endowed with an appropriate $S$--family of data $(D,\omega ,\gamma )$
having the necessary properties  fiberwise.

  \smallskip
 Then we get
 interesting, generally transcendental functions on the base $S$,
 and eventually on moduli spaces/stacks, and these functions satisfy
 (versions of) classical Picard--Fuchs equations.
 
 \medskip
 
 {\bf 1.2.1. Example 1.}  Let $S$ be the affine line with $t$--coordinate,
 and points $t=0, 1$ deleted. Over it, we have the family $\Cal{E}$
 of elliptic curves $E_t$ , that are projective closures of the affine curve
$E_t:\, Y^2=X(X-1)(X-t)$.

 \smallskip
 
 Here is the linear DE for the periods of the relative (over the base)
 1--form
 $dX/Y$  along the closed fiberwise 1--cycles of $E_t$:
  $$
 L_t\omega :=4t(1-t)\frac{d^2\omega}{dt^2} + 4(1-2t)\frac{d\omega}{dt} - \omega =0.
 \eqno(1.2)
 $$

{\bf Example 2.} Non--linear DE  for  the periods of
$dX/Y$ over relative 1--cycles with boundaries
 at sections $P:=(X(t), Y(t))$  of finite order:
 $$
 \mu (P) =0,
 \eqno(1.3)
 $$
where 
$$
\mu (P) := \frac{Y(t)}{2(X(t)-t)^2} - \frac{d}{dt}\left[ 2t(t-1)\,\frac{X^{\prime}(t)}{Y(t)} \right]
+ 2t(t-1)X^{\prime}(t)\,\frac{Y^{\prime}(t)}{Y(t)^2}.
\eqno(1.4)
$$

Notice that $\mu$ defined by $(1.3)$ and extended to the function
on the set of $L$--points of the generic fiber $E_t$ with values in
any differential extension $L$ of $\bold{Q}(t)$
is ``a differential character'': 
$$
\mu (P+Q)=\mu(P) +\mu (Q)
\eqno(1.5)
$$
\smallskip
To explain (and prove) these results, it suffices to notice that
$$
\mu (P)  =L_t \int_{\infty}^P dX/Y
$$
because
$$
L_t (dX/Y) =d \frac{Y}{(X-t)^2}.
$$

\medskip

{\bf 1.3. Perturbative Feynman integrals.} Here I will briefly describe the heuristic origin of
a set of numerical  periods (and periods--functions) indexed  by labeled graphs
relevant for quantum field theory, following [Ma09], sec.~1. For a more focussed study
of (some) of the integrals appearing in this way see [M\"uWZa12] and [W13].

\smallskip

\medskip

 {\it A Feynman path integral} is an heuristic expression of the form
$$
\frac{\int_{\Cal{P}}e^{S(\varphi )}D(\varphi )}{\int_{\Cal{P}}e^{S_0(\varphi )}D(\varphi )}
\eqno(1.6)
$$
or, more generally,  a similar heuristic expression for {\it correlation functions}.

\smallskip

Here the integration domain $\Cal{P}$ stands for a functional space
of {\it classical fields $\varphi$} on  a {\it space--time manifold} $M$. Space--time may be
endowed with a fixed Minkovski or Euclidean metric.
In models of quantum gravity metric is one of the fields. 
Fields  may be scalar functions, tensors of various ranks,
sections of vector bundles, connections.

\smallskip

$S:\,\Cal{P}\to \bold{C}$ is a functional of {\it classical action}: generally $S(\varphi )$
is expressed as an integral over $M$ of a local density on $M$ which is called {\it Lagrangian.}
In our notation (1.6) 
$S(\varphi )=-\int_M L(\varphi (x)) dx.$ Lagrangian density may depend on derivatives, include
distributions etc.  

\smallskip

Usually $S(\varphi )$ is represented as the sum of {\it a quadratic part}
$S_0(\varphi )$ (Lagrangian of free fields) and remaining terms which are 
interpreted as interaction and treated perturbatively.

\smallskip
Finally, the integration measure $D(\varphi )$ and the integral itself $\int_{\Cal{P}}$
should be considered as simply a  part of the total expression (1.6) expressing
the idea of  ``summing the quantum probability amplitudes over all classical trajectories''. 
\smallskip

To explain the appearance and combinatorics of Feynman graphs, we consider a
toy model, in which $\Cal{P}$  is replaced by a finite--dimensional
real space. We endow it with a basis indexed by a finite set
of ``colors''  $A$, and an Euclidean metric $g$ encoded by the symmetric
tensor $(g^{ab}),\,a,b\in A.$ We put $(g^{ab})=(g_{ab})^{-1}.$

\smallskip

The action functional $S(\varphi )$ will be a formal series in linear coordinates on $\Cal{P}$, $(\varphi ^a)$,
of the form
$$
S(\varphi )=S_0 (\varphi) + S_1(\varphi ),\quad
S_0(\varphi ):=-\frac{1}{2} \sum_{a,b} g_{ab}\varphi^a\varphi^b,
$$
$$
S_1(\varphi ):=\sum_{k=1}^{\infty}\frac{1}{k!}\sum_{a_1,\dots ,a_k\in A}
C_{a_1,\dots ,a_k}\varphi^{a_1}\dots \varphi^{a_k}
\eqno(1.7)
$$
where $(C_{a_1,\dots ,a_n})$ are certain  symmetric tensors.
If these tensors vanish for  all sufficiently large ranks $n$,
$S(\varphi )$ becomes a polynomial and can be considered as a genuine function on 
$\Cal{P}$. Below we will treat $(g_{ab})$ and  
$(C_{a_1,\dots ,a_n})$ as independent formal variables, 
``formal coordinates on the space of theories''.

\smallskip

Now we can express  the toy version of (1.6)
as a series over (isomorphism classes of)  graphs.

\smallskip
Here a graph $\tau$  consists of two finite sets, edges $E_{\tau}$
and vertices $V_{\tau}$, and the incidence map   sending $E_{\tau}$
to the set of unordered pairs of vertices. Each vertex is supposed to be incident to at least one edge. There is one {\it empty graph.}

\smallskip

The formula for (1.6) including one more formal parameter $\lambda$ (``Planck's constant'')
looks as follows:
$$
\frac{\int_{\Cal{P}}e^{\lambda^{-1}S(\varphi )}D(\varphi )}{\int_{\Cal{P}}
e^{\lambda^{-1}S_0(\varphi )}D(\varphi )} =
  \sum_{\tau\in\Gamma}\frac{\lambda^{-\chi (\tau )}}{|\roman{Aut}\,\tau |}\,
w(\tau )
\eqno(1.8)
$$
In the right hand side of (1.8), the summation is taken over (representatives of) all 
isomorphism classes
of all finite graphs $\tau$. The weight $w(\tau )$ of such a graph
is determined by the action functional (1.2) as follows:
$$
w(\tau ):=\sum_{u:\,F_{\tau}\to A}\ \prod_{e\in E_{\tau}}
g^{u(\partial e)}\prod_{v\in V_{\tau}} C_{u(F_{\tau}(v))}\,.
\eqno(1.9)
$$
Here  $F_{\tau}$ is the set of  flags, or ``half--edges'' of $\tau$.
Each edge $e$ consists of a pair of flags denoted
$\partial{e}$, and each vertex $v$ determines the set  of flags
incident to it denoted $F_{\tau} (v)$.  Finally, $\chi (\tau )$ is the 
Euler characteristic of $\tau$.

\medskip

The passage of the left hand side of (1.8) to the right  hand side is by definition the result
of term--wise integration of the formal series which can be obtained
from the Taylor series of the exponent in the integrand. Concretely
$$
\int_{\Cal{P}}e^{\lambda^{-1}S(\varphi )}D(\varphi ) =
\int_{\Cal{P}}e^{\lambda^{-1}S_0(\varphi )} \left(1+\sum_{N=1}^{\infty}
\frac{\lambda^{-N}S_1(\varphi )^N}{N!}\right)\, \prod_a d\varphi^a \  :=
$$
$$
\int_{\Cal{P}}e^{\lambda^{-1}S_0(\varphi )} \prod_a d\varphi^a \ +
$$
$$
\sum_{N=1}^{\infty}\frac{\lambda^{-N}}{N!}\sum_{k_1,\dots ,k_N=1}^{\infty}\frac{1}{k_1!\dots k_N!}
\sum_{a^{(i)}_j\in A, 1\le j\le k_i}\prod_{i=1}^N
C_{a_1^{(i)},\dots ,a_{k_i}^{(i)}}
\int_{\Cal{P}} e^{\lambda^{-1}S_0(\varphi )}\prod_{i,j}^N \varphi^{a_j^{(i)}}
 \prod_a d\varphi^a  \,.
  \eqno(1.10)
$$
\smallskip

This definition makes sense if the right hand side of (1.10) is understood
as a formal series of infinitely many independent weighted variables $C_{a_1,...,a_k}$,
weight of  $C_{a_1,...,a_k}$ being $k$. In fact, the Gaussian
integrals in the coefficients uniformly converge, and one can use the so called
Wick's lemma.
\smallskip

The last remark is that {\it periods} appearing in concrete models of quantum field theories
are {\it weights (1.9)}, in which the {\it summation} over maps $u:\,F_{\tau}\to A$
is replaced by the {\it integration} over some continuous variables such as positions/momenta/colours 
of  particles moving along the edges of the respective Feymann graph: cf. [W13],
[M\"uWZa12] and references therein.

\bigskip


\centerline{ \bf 2. Arithmetic differential equations}

\medskip

{\bf 2.1. Analogies between $p$--adic numbers and formal series.} Combining the lessons of previous examples we suggest now that
in order to see  ``$p$--adic properties'' of numerical periods,
transcendental numbers important for physics, one could try to design
a theory of  ``derivations in $p$--adic direction'' and interpret  numerical periods
as solutions of  differential equations in  the $p$--adic direction.

\smallskip

Below we present basics of such a theory due to A.~Buium. We start with the following
table of analogies. On the formal series side, we consider rings of the form $k[[t]]$
where $k$ is a field of characteristics zero. On the $p$--adic side, we consider the
maximal unramified extension $R$ of $\bold{Z}_p$.


\newpage



\quad \quad $\underline{POWER\ SERIES}$\quad\quad\quad\quad\quad\quad \quad\quad\quad $\underline{p-ADICS}$

\bigskip

\quad\quad\quad$ \sum a_it^i\in k[[t]]=:L$  \quad\quad\quad\quad\quad\quad\quad\quad   $ \sum \varepsilon_i p^i\in R := \bold{Z}_p^{un}$

\bigskip

\quad Field of constants: $a_i\in k$ \quad\quad\quad\quad Monoid: $\varepsilon_i\in \mu_{\infty}\cup \{0\}$

\smallskip

\quad\quad\quad\quad\quad \quad\quad\quad\quad\quad\quad\quad\quad\quad\quad {\it (Teichm\"uller representatives)}

\bigskip

\quad Derivation: $d/dt$  \quad\quad\quad\quad\quad\quad\quad\quad\quad $\delta_p(*):=\frac{\Phi(*)-*^p}{p}$
\smallskip

\quad\quad\quad\quad\quad \quad\quad\quad\quad\quad\quad\quad\quad \quad\quad\quad $(\Phi:= lift\ of\ Frobenius)$
\bigskip
Polynomial Diff. Operators (PDO): \quad\quad \quad\quad  $p$-adic PDO:

\medskip

\quad$D\in L[T_0,T_1,\dots ,T_n]$ \quad\quad\quad\quad\quad\quad  $D_p\in \overline{R[T_0,T_1,\dots ,T_n]}$
\smallskip

\quad\quad\quad\quad\quad\quad\quad\quad\quad\quad\quad\quad\quad\quad\quad\quad\quad {\it ($p$--adic completion!)}

\bigskip

\centerline{----------------------------------------------------------------------------}
\centerline{Action of PDO: $f\mapsto D(f, f^{\prime},\dots f^{(n)})$ or $D_p(f,\delta_pf,\dots ,\delta^n_pf)$}
\centerline{----------------------------------------------------------------------------}

\bigskip

The Frobenius lift $\Phi :\,R\to R$ involved in the definition of the $p$--adic derivative $\delta_p$
is given explicitly as $\Phi (\sum \varepsilon_i p^i) := \sum \varepsilon_i^p\, p^i$.

\bigskip

{\bf 2.2. Examples and applications.} Here we  give a sample of interesting $p$--adic differential
operators. 

\medskip

{\bf 2.2.1. Example 1: {\it $p$--adic logarithmic derivative.}} It is an analog of  the map
$$
\bold{G}_m(L) \to\bold{G}_a(L):\quad f\mapsto f^{\prime}/f
\eqno(2.1)
$$
where a point $x\in \bold{G}_m(L)$ is represented by the value $f\in L^*$ at $x$ of a fixed
algebraic character $t$ of $\bold{G_m}$ such that $\bold{G_m}=\roman{Spec}\,[t,t^{-1}]$.
Similarly, its $p$--adic version is the differential character 
$$
\bold{G}_m(R) \to\bold{G}_a (R):
$$
$$
\quad a\mapsto \delta_pa\cdot a^{-p}-\frac{p}{2}(\delta_pa\cdot a^{-p})^2
+\frac{p^2}{3}(\delta_pa\cdot a^{-p})^3 - \dots
\eqno(2.2)
$$

\medskip

{\bf Example 2: {\it Quadratic reciprocity symbol:}}
$$
\left(\frac{a}{p}\right) = a^{\frac{p-1}{2}}\left(1+\sum_{k=1}^{\infty}
(-1)^{k-1}\frac{(2k-2)!}{2^{2k-1}(k-1)! k!}  (\delta_pa)^k a^{-pk}\right) \, .
$$

\medskip

{\bf Example 3:} {\it a $p$--adic analog of the differential character $\mu$ of the group
of sections of a generic elliptic curve:}
$$
\mu(P) = 
(4t(1-t)\frac{d^2}{dt^2} + 4(1-2t)\frac{d}{dt} - 1)\int_{\infty}^P \frac{dX}{Y}
 $$
as a non--linear $p$--adic  DO acting upon coordinates of $P$.

\smallskip

Such analogs were constructed in [Bu95] also for abelian varieties of arbitrary dimension
and called $\delta_p$--differential characters $\psi ( P )$. More precisely, let $E$ be an elliptic
curve over $R$.  Then there exist a differential additive map $\psi :\, E( R )\to R^+$ of order 2 (as 
in the geometric case) or 1 (as for $\bold{G}_m$). 
\smallskip
A character of order
2 exists if $E$ has a good reduction and is not the canonical lift of its reduction in the sense of Serre--Tate:
cf. additional discussion in 4.4 below.

\smallskip

A character of order 1 exists if either $E$ has good ordinary reduction and is the canonical lift,
or $E$ has a bad multiplicative reduction.

\smallskip

Using these multiplicative characters, A.~Buium and  the author constructed in [BuMa13]
{\it ``Painlev\'e VI equations with $p$--adic time.''}

\medskip
{\bf 2.3. General formalism of $p$--derivations.} In the commutative algebra, given a ring $A$ and an $A$--module $N$,
{\it a derivation of $A$ with values in  $N$ } is any additive map 
$\partial :\,A\to N$ such  that $\partial (ab)= b\partial a+ a \partial b$. Equivalently,
the map $A\to A\times N:\, a\mapsto (a,\partial a)$
is a ring homomorphism, where $A\times N$ is endowed with the structure
of commutative ring with componentwise addition, inheriting multiplication from $A$ on $A\times \{0\}$
and having $\{0\}\times N$ as an ideal of square zero.

\smallskip

Similarly, in arithmetic geometry Buium defines  {\it a $p$--derivation of $A$ with values in an $A$--algebra $B$,
$f:A\to B$, }
as a map $\delta_p:\,A\to B$ such that the map $A\to B\times B:\, a\mapsto (f(a),\delta_p( a))$
is a ring homomorphism $A\to W_2(B)$ where $W_2(B)$ is the ring of $p$--typical Witt vectors of length 2.
Here Witt vectors of the form $(0,b)$ form the ideal of square zero only if  $pB=\{0\}$.
 
 \smallskip
Making this definition explicit, we get $\delta_p(1)=0$, and the following versions of additivety
and Leibniz's formula:
$$
\delta_p(x+y)=\delta_p(x)+\delta_p(y) +C_p(x,y),
\eqno(2.3)
$$
$$
\delta_p(xy)=f(x)^p\cdot \delta_p(y)+f(y)^p\cdot \delta_p(x) +p\cdot \delta_p(x)\cdot\delta_p(y),
\eqno(2.4)
$$
where
$$
C_p(X,Y):=\frac{X^p+Y^p-(X+Y)^p}{p}\in \bold{Z}[X,Y].
\eqno(2.5)
$$
In particular, this implies that for any $p$--derivation $\delta_p:\,A\to B$ the respective map $\phi_p :\, A\to B$
defined by $\phi_p (a):=f(a)^p+p\delta_p(a)$ is a ring homomorphism
satisfying $\phi_p (x)\equiv f(x)^p\,\roman{mod}\,p$, that is ``a lift of the Frobenius map applied to $f$''.
\smallskip

Conversely, having such a lift of Frobenius, we can uniquely reconstruct the respective derivation $\delta_p$
{\it under the condition that $B$ has no $p$--torsion:}
$$
\delta_p(a):= \frac{\phi_p(a) - f(a)^p}{p}
$$
generalising the  definition given in 2.1 for $A=B=R$ and identical morphism.

\smallskip

Working with $p$--derivations $A\to A$ with respect to the identity map $A\to A$
and keeping $p$ fixed, we may call $(A,\delta )$ 
 a $\delta$--ring. Morphisms of $\delta$--rings are algebra morphisms
 compatible with their $p$--derivations.

\medskip

{\bf 2.4. $p$--jet spaces.}  Let $A$ be an $R$--algebra. {\it A prolongation sequence} for $A$
consists of a family of $p$--adically complete $R$--algebras $A^i, i\ge 0$, 
where $A^0=A\,\widehat{}$\quad  is the $p$--adic comple\-tion of $A$, and of maps 
$\varphi_i, \delta_i:\,A^i\to A^{i+1}$ satisfying the following conditions:

\smallskip

{\it a) $\varphi_i$ are ring homomorphisms, each $\delta_i$ is a $p$--derivation  with
respect to $\varphi_i$, compatible with $\delta$ on $R$.
\smallskip

b) $\delta_i\circ \varphi_{i-1}=\varphi_i\circ \delta_{i-1}$ for all $i\ge 1$.}
\smallskip

Prolongation sequences  form a category with evident morphisms, ring homomorphisms
$f_i:A^i\to B^i$ commuting with $\varphi_i$ and $\delta_i$,
and in its subcategory with fixed $A^0$ {\it there exists an initial element,
defined up to unique isomorphism} (cf. [Bu05], Chapter 3). It can be called the
universal prolongation sequence.

\smallskip

In the geometric language, if $X=\roman{Spec}\,A$, the formal spectrum of the $i$--th ring $A^i$
in the universal prolongation sequence  is denoted $J^i(X)$ and called 
{\it the $i$--th $p$--jet space of $X$.} Conversely, $A^i=\Cal{O}(J^i(X))$, the ring of global
functions. 

\smallskip
The geometric morphisms (of  formal schemes over $\bold{Z}$) corresponding to $\phi_i$ are denoted
$\phi^i:\, J^i(X)\to J^0(X)=:X\,\widehat{}$ (formal $p$--adic completion of $X$).

\smallskip

This construction is compatible with localisation so that
it can be applied to the non--necessarily affine schemes: cf. [Bu05], Chapter 3.

\bigskip

\centerline{\bf 3. An arithmetically global version of Buium's calculus}

\smallskip

\centerline{\bf and lambda--rings}

\medskip
 {\bf 3.1. Introduction.}  $p$--adic numbers  were considered in sec. 2 above as analogs of formal functions/local germs of functions
of one variable.
\smallskip
In this section, we discuss the following question: does there exist a (more) global version of ``arithmetic functions'', 
elements of a ring $A$, admitting $p$--adic derivations $\delta_p$ with respect to 
 several, eventually all primes $p$?

\medskip

An obvious example is  $\bold{Z}$:
$$
\delta_p(m)= \frac{m-m^p}{p}.
$$

Generally, we need ``lifts of Frobenii'': such ring endomorphisms $\Phi_p:\,A\to A$
that $\Phi (a) \equiv a^p\,\roman{mod}\,p$. Then we may put
$$
\delta_p(a)= \frac{\Phi_p(a)-a^p}{p}.
$$

A general framework for a coherent system of such lifts is given by the following definition:

\medskip

{\bf 3.2. Definition.}  {\it A system of  psi--operations on a commutative
unitary ring $A$ is a family of {\it ring endomorphisms} $\psi^k:\,A\to A$, $k\ge 1$, such that:
$$
\psi^1=id_A,\quad \psi^k\psi^r=\psi^{kr},
$$
$$
\psi^px\equiv x^p\,\roman{mod}\,pA\quad for\ all\ primes\ p.
$$}

Another important structure is introduced by the following definition:
\medskip
{\bf 3.3. Definition.} {\it A system of {\it lambda--operations} on a commutative
unitary ring $A$ is a family of {\it additive group endomorphisms} $\lambda^k:\,A\to A$, $k\ge 0$,
such that
$$
\lambda^0(x)=1,\ \lambda^1=id_A,
$$
$$
\lambda^n(x+y)=\sum_{i+j=n} \lambda^i(x)\lambda^j(y).
$$
}
  These structures are related in the following way:
  
  \medskip

  {\bf 3.4. Proposition.}  {\it (a) If $A$ has no additive torsion, then any system of psi--operations defines a unique system of
lambda--operations satisfying  the compatibility relations:
$$
(-1)^{k+1}k\lambda^k(x)=\sum_{i+j=k,\, j\ge 1}(-1)^{j+1}\lambda^i(x)\psi^j(x) .
$$

(b) Generally, any system of lambda--operations defines a unique system of
psi--operations satisfying the same compatibility relations.}

\bigskip

Briefly, such a ring, together with psi's and lambda's, is called {\it a lambda--ring.}

\medskip
{\bf 3.5. Example: a Grothendieck ring.} 
Let $\Cal{R}$ = a commutative unitary ring.

\smallskip

Denote by  $A=A_{\Cal{R}}$  the Grothendieck $K_0$--group of the additive category, consisting of pairs
$(P, \varphi )$, where $P$ is a projective $\Cal{R}$--module of finite type, $\varphi :\,P\to P$
an endomorphism. Denote by $[(P,\varphi )]\in A$ the class of $(P,\varphi ).$

\smallskip
The ring structure on $A$ is induced by the tensor product: $[(P,\varphi )][(Q,\psi )] :=  [(P\otimes Q,\varphi \otimes \psi ]$.

\smallskip

The lambda--operations on $A$ are defined by $\lambda^k\,[(P,\varphi )]:=  [(\Lambda^kP,\Lambda^n \varphi )]$.

\medskip

{\bf 3.6. Example: the big Witt ring $W(\Cal{R})$.} 
Again, let $\Cal{R}$ = a commutative unitary ring.

\smallskip

Define the additive group of $W(\Cal{R})$ as {\it the multiplicative group} $1+T\Cal{R}[[T]]$.

\smallskip

The multiplication $*$ in $W(\Cal{R})$ is defined on elements $(1-at)$ as $(1-aT)*(1-bT):= 1-abT$, and then extended 
to the whole $W(\Cal{R})$ by distributivity, continuity in the $(T)$--adic topology, and functoriality in $\Cal{R}$.

\medskip

Similarly, lambda--operations  in $W(\Cal{R})$ are defined by $\lambda^k\,(1-aT):= 0$ for $k\ge 2$, and then
extended by addition formulas (Def. 3.3)  and  continuity.

\bigskip

\centerline{\bf 4. Roots of unity as constants:}

\medskip

\centerline{\bf geometries over  ``fields of characteristic 1''}

\bigskip
{\bf 4.1. Early history.} In the paper [T57], J.~Tits noticed that some basic numerical invariants
related to the geometry of classical groups over finite fields $\bold{F}_q$ have  well--defined
values for $q=1$, and these values admit suggestive combinatorial interpretations.
\smallskip

For example, if $q=p^k$, $p$ a prime, $k\ge 1$, then
$$
\roman{card}\,\bold{P}^{n-1}(\bold{F}_q)=
 \dfrac{\roman{card}\,(\bold{A}^n(\bold{F}_q)\setminus \{0\})}{\roman{card}\,\bold{G}_m(\bold{F}_q)}=\dfrac{q^n-1}{q-1} =: [n]_q,
 $$
$$
\roman{card}\,Gr\,(n,j)(\bold{F}_q) = \roman{card}\,\{
\bold{P}^j(\bold{F}_q) \subset\bold{P}^n(\bold{F}_q)\}=: {\binom nj}_q,
$$
and the $q=1$ values of the right hand sides are cardinalities of the sets
\medskip
\centerline{$\bold{P}^{n-1}(\bold{F}_1)$:= a finite set $P$ of cardinality $n$,}
\smallskip
\centerline{$Gr\,(n,j)(\bold{F}_1)$ := the set of subsets of $P$ of cardinality $j$.}
\medskip

Tits suggested a program: make sense of algebraic geometry over 
{\it ``a field of characteristic one''}
so that the ``projective geometry'' above becomes a special case
of the geometry of Chevalley groups and their homogeneous spaces.

\smallskip
The first implementation of Tits' program was achieved only
in 2008 by  A.~Connes and  C.~Consani, cf. [CC11], after the foundational
work by C.~Soul\'e  [So04]. However, they required $\bold{F}_{1^2}$ as a definition field.
\smallskip

Earlier, in an unpublished manuscript [KaS], M.~Kapranov and A.~Smirnov introduced
fields $\bold{F}_{1^n}$ on their own right.

\smallskip

They defined $\bold{F}_{1^n}$ as the monoid $\{0\}\cup \mu_n$, where $\mu_n$
is the set of roots of unity of order $n$. Moreover, they defined a
a vector space over $\bold{F}_{1^n}$ as a pointed set $(V,0)$
with an action of $\mu_n$ free on $V\setminus \{0\}$.
The group $GL(V)$, by definition, consists of permutations of $V$ compatible with action of $\mu_n$.
Kapranov and Smirnov defined the determinant map  $\roman{det}:\, GL(V)\to \mu_n$
and proved a beautiful formula for the power residue symbol.
\smallskip

Namely, if $q=p^k\equiv 1\,\roman{mod}\,n$ and $\mu_n$ is embedded
in $\bold{F}_q^*$, $\bold{F}_q$ becomes a vector space over
$\bold{F}_{1^n}$, and the power residue symbol
$$
\left(\dfrac{a}{\bold{F}_q}\right)_n:= a^{\frac{q-1}{n}}\in \mu_n
$$
is {\it the determinant of the multiplication by $a$} in   $\bold{F}_{1^n}$--geometry.
\smallskip
 
 Cf. also [Sm92], [Sm94].

\medskip
As we noticed in sec. 2, constants with respect to Buium's derivation 
$\delta_p$ in $R:= \bold{Z}_p^{un}$ are roots of unity (of degree prime to $p$) completed by $0.$ 
\smallskip 
Therefore, in the context of the differential geometry ``in the $p$--adic direction''  
an independent project of
Algebraic Geometry
``over roots of unity'', or ``in characteristic 1'', or else
``over fields $\bold{F}_1, \bold{F}_{1^n},  \bold{F}_{1^\infty}$'' É
acquires a new motivation. Moreover, it becomes enriched with new insights: whereas
at the first stage schemes in characteristic 1 were constructed by glueing ``spectra of commutative monoids'',
now they could be conceived as $\bold{Z}$--schemes endowed with lambda--structure
considered as descent data: see [Bor11a,b], [Bor09]. Here is a brief survey of Borger's philosophy,
showing that his schemes form a natural habitat for $p$--adic differential geometries
as well.
\medskip

{\bf 4.2. Borger's philosophy.} The category of affine $\bold{F}_1$--schemes $Aff_1$
can be defined as the opposite  category  of rings endowed with
lambda--structures, $(A,\Lambda_A)$, and compatible morphisms. The forgetful functor
to the usual category of affine schemes  $Aff_1\to Aff:\, (A,\Lambda)\mapsto A$
is interpreted as the functor of base extension $*\mapsto *\otimes_{\bold{F}_1}\bold{Z}$.
\smallskip

Thus, a lambda--structure on a ring $A$ is {\it a descent data on} $\roman{Spec}\, A$ to $\bold{F}_1.$

\smallskip

 In particular, $W(\bold{Z})$ must be considered as (a completion of?) 
 $\bold{Z}\otimes_{\bold{F}_1} \bold{Z}$.
 
\smallskip

More\ generally, using general topos theory, Borger globalizes this construction, constructing a natural {\it algebraic geometry of $\lambda$--schemes,} which should be thought of as a lifted {\it algebraic geometry over $\bold{F}_1.$} 

\smallskip

Just as all of usual algebraic geometry is contained in the big \'etale topos of $\bold{Z}$, $\lambda$--algebraic geometry is contained in a big topos, which should be thought of as the big \'etale topos over $\bold{F}_1$.  There is a map of topoi from the big etale topos over $\bold{Z}$ to the one over 
 $\bold{F}_1.$
 \smallskip
Schemes of finite type over $\bold{F}_1$ (in this sense, as in most other approaches) are very {\it rigid, combinatorial objects.}  They are essentially quotients of toric varieties by toric equivalence relations.  

\smallskip
Non--finite--type schemes over 
 $\bold{F}_1$ are  more interesting.  The big de Rham--Witt cohomology of $X$ 
 ``is'' the de Rham cohomology of $X$ ``viewed as an $\bold{F}_1$--scheme". It should contain the full information of the motive of $X$ and is probably a concrete universal Weil cohomology theory.  
 
 \smallskip
 
 The Weil restriction of scalars from $\bold{Z}$ to $\bold{F}_1$ exists and is an arithmetically global version of Buium's $p$--jet space.
\smallskip

In conclusion, we briefly mention some remaining challenges. 
 
 \medskip
 
 {\bf 4.3. Euler factors at infinity and $\bold{F}_1$--geometry.} In [Ma95], I suggested that
 there should exist a category of $\bold{F}_1$--motives visible through the $q=1$ point count of
 $\bold{F}_1$--schemes. Predictions about such a point count were justified in
 Soul\'e's geometry, cf. [So04].   In particular
 the zetas of
non--negative  powers of  the ``Lefschetz (dual Tate) motive'' 
$\bold{L}$ must be:
$$
Z(\bold{L}^{\times n},s)= \frac{s+n}{2\pi}.
$$

This provides a conjectural {\it bridge between $\bold{F}_1$--geometry}
and geometry of $Spec\,\bold{Z}$ at the archimedean infinity,
that is, {\it Arakelov geometry}:
a $\Gamma$--factor of  classical zetas, e.g.,
$$
\Gamma_\bold{C}(s):=[(2\pi )^{-s}\Gamma (s)]^{-1}=\prod_{n\ge 0} \frac{s+n}{2\pi}
$$
(regularized product)
looks like $\bold{F}_1$--zeta of the dualized  inf--dim  projective space over $\bold{F}_1$.
\smallskip

However, this phenomenon  remains an isolated  observation, and 
the archimedean prime  still remains ``first among equals'' breaking
the democracy of all valuations.

\medskip

{\bf 4.4. Other geometries} ``under $\roman{Spec}\,\bold{Z}$''.  In the traditional algebraic geometry,
the special role of $\roman{Spec}\,\bold{Z}$ is related to the fact that it is the final object
of the category of schemes.  Since it is very far from being ``a point--like object'',
it seemed natural to imagine that $\roman{Spec}\,\bold{F}_1$, being ``really point--like'',
will replace it. However,  the belief that in an extended algebraic geometry
there should necessarily exist a final object, is  unfounded. Already
in the simplest category of Deligne--Mumford stacks over a field $k$, admitting quotients with respect
to the trivial action of any finite group $G$, there is no final object, because we have non--trivial morphisms 
$\roman{Spec}\,k\to\roman{Spec}\,k/G$.

\smallskip

This led several authors to the contemplation of more general geometries
lying  ``under $\roman{Spec}\,\bold{Z}$'' but not necessarily at the  bottom of the unfathomable abyss:
cf. the To\"en--Vaqui\'e project  [TV05].

\smallskip

For example, in the Borger--Buium's framework we may consider
schemes for which Frobenius lifts are given only for some subsets of primes,
eventually one prime $p$, such as the Serre--Tate canonical liftings
of Abelian varieties in characteristic $p$: cf. [Katz81].

\smallskip

More precisely,  for the simplest case of elliptic curves,   
denote by $M$  the $p$--adic completion of the moduli stack of elliptic curves 
without  supersingular locus.  One can define Frobenius lift on this stack: it sends an elliptic curve to its quotient by its canonical subgroup. The latter is defined as the unique closed sub--groupscheme whose Cartier dual is the \'etale lift to 
$\bold{Z}_p$ of the Cartier dual of the kernel of Frobenius on the fiber modulo $p$.   This endomorphism also lifts to a natural endomorphism of the universal elliptic curve. So James Borger suggests to say that $M$ ``descends to the $p$--typical $F_1$'', and the same can be said about the universal elliptic curve over it. The $p$--adic elliptic curves
with Frobenius lift are called canonical liftings.

\smallskip

Notice that if we replace the $p$--adic direction by the functional one, we would simply speak
about families of elliptic curves with constant absolute invariants. But $p$--adic absolute invariants
of canonical liftings are by no means ``constants'' in the naive sense, discussed in sec. 2, that is
they are not Teichm\"uller  representatives: cf. a recent paper by Finotti,
"Coordinates of the $j$--invariant of the canonical lifting", posted at
 http://www.math.utk.edu/~finotti/ , and [Er13].

\smallskip

A better understanding of this discrepancy presents
 an interesting challenge for the $p$--adic differential geometry.
 
 \smallskip
 
 \smallskip

{\it Acknowledgements.} Collaboration with A.~Buium on [BuMa13] helped me much in conceiving
this survey. J.~Borger generously explained me some of his constructions and motivations.
Igor Volovich stimulated the final writing by inviting me to give a talk at
the  International Workshop on $p$--adic methods for modelling of
complex systems, Bielefeld, April 15--19, 2013. I am grateful to them all.

\vskip1cm

\centerline{\bf References}

\medskip

[ACG13] A.~Abdessalam, A.~Chandra, G.~Guadagni. {\it Rigorous quantum field theory
functional integrals over the $p$--adics I: anomalous dimensions.} arXiv:1302.5971

\smallskip

[Be11] A.~Beilinson. {\it $p$--adic periods and derived De Rham cohomology.}
Journ. AMS, vol. 25, no. 3 (2012), 319--327. arXiv:1102.1294

\smallskip

[Bor09] J.~Borger. {\it Lambda--rings and the field with one element.} arXiv:0906.3146
\smallskip

[BorBu09] J.~Borger, A.~Buium.
{\it Differential forms on arithmetic jet spaces.} Selecta Math. (N.S.) 17 (2011), no. 2, 301--335.
arXiv:0908.2512

[Bor11a] J.~Borger. {\it The basic geometry of Witt vectors, I: The affine case.}
 Algebra Number Theory 5 (2011), no. 2, 231--285.

\smallskip
[Bor11b] J.~Borger. {\it  basic geometry of Witt vectors. II: Spaces.} Math. Ann. 351 (2011), no. 4, 877--933.
\smallskip

[Bu95] A.~Buium. {\it Differential characters of Abelian varieties over p--adic fields.}
Inv.~Math., vol. 122 (1995), 309--340.

\smallskip
[Bu05] A.~Buium. {\it Arithmetic Differential Equations.} AMS Math Surveys
and Monographs, vol. 118, 2005.
\smallskip

[BuMa13] A.~Buium, Y.~Manin. {\it Arithmetic Differential Equations of Painlev\'e  VI Type.}
arXiv:1307.3841

\smallskip

[CC11] A.~Connes, C.~Consani. {\it On the notion of geometry over $\bold{F}_1$}. J. Algebraic Geom. 20 (2011), no. 3, 525--557.
\smallskip

[CCMa08] A.~Connes, C.~Consani, M.~Marcolli. {\it Fun with $F_1$.}
J. Number Theory 129 (2009), no. 6, 1532--1561.  math.AG/0806.2401

\smallskip

[De04] A.~Deitmar. {\it Schemes over $F_1$.} In: Number Fields and Function Fields -- Two Parallel
Worlds. Ed. by G.~ van der Geer, B.~Moonen, R.~Schoof. Progr. in Math, vol. 239, 2005.  math.NT/0404185

\smallskip

[Er13] A.~Erdogan. {\it A universal formula for the $j$--invariant of the canonical
lifting.} arXiv:1211.1152

\smallskip

[GoMa04] A.~Goncharov, Yu.~Manin. {\it Multiple zeta--motives and moduli spaces $\overline{M}_{0,n}$.}
Compos. Math. 140:1 (2004), 1--14. math.AG/0204102

\smallskip

[Fa88] G.~Faltings. {\it $p$--adic Hodge theory.} J.~Amer.~Math.~Soc., 1(1988), 255--288.

\smallskip
[KaS] M.~Kapranov, A.~Smirnov. {\it Cohomology determinants and reciprocity laws:
number field case.} Unpublished manuscript, 15 pp.
\smallskip

[Katz81] Katz, N. {\it Serre--Tate local moduli.} Algebraic surfaces (Orsay, 1976--78), pp. 138--202, Lecture Notes in Math., 868, Springer, Berlin-New York, 1981.

\smallskip

[KoZa01] M.~Kontsevich, D.~Zagier. {\it Periods.} In: Mathematics unlimited---2001
and beyond, 771--808, Springer, Berlin, 2001.
\smallskip

[LeBr13] L.~Le Bruyn. {\it Absolute geometry and the Habiro topology.}
arXiv:1304.6532

\smallskip

[Ma89] Yu.~Manin. {\it Reflections on arithmetical physics.} In: Conformal
Invariance and string theory (Poiana Brasov, 1987),
Academic Press, Boston, MA, 1989, 293--303. Reprinted in
``Mathematics as Metaphor'', Selected Essays by Yu.~I.~Manin, AMS 2007,  pp. 149--155.

\smallskip
[Ma95] Yu.~Manin. {\it Lectures on zeta functions and motives (according to Deninger and Kurokawa).}
Ast\'erisque 228:4 (1995), 121--163.

\smallskip

[Ma08] Yu.~Manin. {\it Cyclotomy and analytic geometry over $\bold{F}_1$.} In: Quanta of Maths. Conference in honour of Alain Connes.
Clay Math. Proceedings, vol. 11 (2010), 385--408.
Preprint math.AG/0809.2716.

\smallskip

[Ma09] Yu.~Manin. {\it Renormalization and computation I: motivation and background.}
In: Proceedings OPERADS 2009, eds. J. Loday and B. Vallette,
S\'eminaires et Congr\`es 26, Soc. Math. de France, 2012, pp. 181--223.
Preprint math.QA/0904.4921

\smallskip

[M\"uWZa12] S.~M\"uller--Stach, S.~Weinzierl, R.~Zayadeh. {\it Picard--Fuchs equations for
Feynman integrals.} arXiv:1212.4389

\smallskip
[Sm92] A.~L.~Smirnov. {\it Hurwitz inequalities for number fields. (Russian).} Algebra i Analiz 4 (1992), no. 2, 186--209;
translation in St. Petersburg Math. J. 4 (1993), no. 2, 357--375.

\smallskip

[Sm94] A.~L.~Smirnov. {\it Absolute determinants and Hilbert symbols.} Preprint MPI 94/72, Bonn, 1994.
\smallskip
[So04] C.~Soul\'e. {\it Les vari\'et\'es sur le corps \`a un \'el\'ement.} Mosc. Math. J. 4:1 (2004), 217--244.
\smallskip

[Ti57] J.~Tits. {\it Sur les analogues alg\'ebriques des groupes semi--simples complexes.}
Colloque d'alg\`ebre sup\'erieure, Centre Belge de Recherches Math\'ematiques,
\'Etablissement Ceuterick, Louvain, 1957, 261--289.

\smallskip

[TV05] B.~To\"en, M.~Vaqui\'e. {\it Au--dessous de $Spec\,\bold{Z}$.} J. K-Theory 3 (2009), no. 3, 437--500.  math.AG/0509684

\smallskip

[VlVoZe94] V.~S.~Vladimirov, I.~V.~Volovich, E.~I.~Zelenov. {\it p--adic analysis and mathematical physics.}
Series on Soviet and East European Math., 1. World Scientific, River Edge, NJ, 1994.

\smallskip

[W13] S.~Weinzierl. {\it Periods and Hodge structures in perturbative quantum field theory.}
arXiv:1302.0670 [hep--th]

\enddocument